*Article*

# Enhancing Complex Injection Mold Design Validation Using Multicombined RV Environments

Jorge Manuel Mercado-Colmenero [1,2], Diego Francisco Garcia-Molina [1], Bartolomé Gutierrez-Jiménez [3] and Cristina Martin-Doñate [1,2,*]

1   Department of Engineering Graphics Design and Projects, University of Jaen, Campus Las Lagunillas s/n, Building A3-210, 23071 Jaén, Spain; jmercado@ujaen.es (J.M.M.-C.); dfgarcia@ujaen.es (D.F.G.-M.)
2   INGDISIG Jaén Research Group, Campus Las Lagunillas s/n, Building A3, 23071 Jaén, Spain
3   Moldes Tuccibérica S.L., C/Matriceros s/n Aptd nº 88, 23600 Jaén, Spain; moldestucciberica@gmail.com
*   Correspondence: cdonate@ujaen.es; Tel.:+34-953-212-821; Fax: +34-953-212-334

**Abstract:** The intricate design of real complex injection molds poses significant challenges. Mold design validation often falls to operators with tool-handling experience but limited CAD proficiency. Unlike other industries, the scale and costs of injection mold fabrication hinder prototyping before production. Virtual reality (VR) has emerged as a revolutionary solution offering a safe, immersive, and realistic experience and accessible using QR codes. This paper presents a new multimodal virtual environment tailored to validate mold design complexities. Integrating knowledge-enriched visual tools like interactive 3D models and dynamic visualizations enables users to explore complex mold designs. Statistical analyses, including the Wilcoxon test, unveil significant differences in interference detection, internal topology tracking, and validation of assembly and disassembly accessibility for both small and large mold components when comparing validation conducted through traditional means using solely CAD systems versus the utilization of multidimensional validation methods. Efficiency gains in using VR devices for mold design validation in a hybrid environment in the analysis of relative frequencies. The present study surpasses the state of the art illustrating how VR technology can substantially reduce manufacturing errors in injection molding processes, thereby offering important advantages for manufacturers emerging as an essential tool for this impact industry in the next years.

## 1. Introduction

Plastic injection molding is currently the most widespread plastic transformation process globally, being crucial in the manufacturing of a wide range of consumer products. It is estimated that more than 70% of these products include components produced using this technology [1]. In 2020, the global injection molded plastic market reached USD 265.1 billion, with a forecast annual growth of 4.6% until 2028 [2]. This production process is distinguished by its rigorous specifications, demanding resistant and long-lasting tooling to guarantee the production of high-quality parts in an optimal time. The geometric complexity of plastic parts manufactured using injection molding, with free shapes, clips, bosses and stiffeners, implies challenges in demolding, requiring sophisticated mechanical systems composed of multiple pieces of various materials that must function in a synchronized manner [3]. Furthermore, given the high pressure during the process, molds usually divide the cavity into interchangeable parts to avoid deformations in small part details due to bending stresses that would require changing the entire cavity [4]. In the field of plastic part manufacturing, molds play a crucial role, comprising various fixed and movable components [5]. This complexity demands a thorough understanding and detailed examination of both internal and external surfaces within the assembly, which poses challenges in design and review validation. Additionally, manufacturing injection molds involves significant costs and lengthy production times, which limit prototyping options before mass production [6]. Validating design injection molds is a critical process that ensures the accuracy and functionality of the mold before widespread use [7]. The globalization of injection mold manufacturing further complicates this issue, as molds are often designed in different countries, requiring comprehensive design reviews to prevent communication breakdowns during production. Furthermore, ensuring there is no interference among mold components is crucial to prevent damage during manufacturing and maintain the quality of the parts [8–10]. Validation of assembly and disassembly, alongside component arrangement within the mold, is essential to the design and manufacturing process [11]. A poorly designed or difficult-to-disassemble mold can increase injection machine downtime and impede production line efficiency. Ensuring accessible mold components for inspection, maintenance, and repair is vital, especially for critical elements like ejection mechanisms. Improper assembly or layout can increase the risk of premature mold damage or wear, leading to costly downtime and a shortened mold life. Thus, prioritizing accessibility and implementing efficient maintenance procedures are crucial for optimizing mold performance and minimizing production disruptions [12].

Virtual reality (VR) has emerged as a promising tool for revolutionizing the validation process of injection molds [13]. By creating immersive virtual environments that accurately replicate mold operations, VR enables users to visualize and analyze internal components with exceptional precision [14]. It has therefore been decided to share the experience through dissemination channels with access via QR codes. This resource substantially improves the interpretation of the CAD model and its operation. Therefore, different operators will be able to access CAD modelling in an immersive and intuitive way through the aforementioned VR resource for the compression of the mold, identification of problems and facilitation of their information by performing CAD modelling more easily. Moreover, VR provides a safe testing environment, resulting in significant cost and time savings compared to traditional methods [15]. In contrast to the limitations of computer-aided design (CAD) programs, VR offers a more intuitive and dynamic experience. While 2D and 3D CAD on screens have been widely used, they often fall short in functional and ergonomic validation [16]. CAD software lacks the intuitive interface necessary for users without CAD or computer science experience [17–19]. Conversely, VR introduces innovative possibilities for data interaction, enabling engineers to visualize projects in 3D and gain deep insights into their functionality. Interactive assembly and disassembly simulations facilitated by VR simplify the comprehension of complex processes [20]. This capability allows engineers to engage with their designs more effectively, leading to an improved understanding and optimization of mold functionalities.

For highly complex injection molds, particularly those with over 300 components, depending on interactivity within VR environments may prove inadequate. Exploring additional methods is essential when addressing mold designs, as VR, while promising, may have limitations, especially with complex designs. Therefore, designing multimodal environments tailored to mold complexities is crucial for achieving effective outcomes. One approach involves creating multiple scenarios analyzing various design and functionality aspects. Integrating knowledge-enriched visual tools is vital for deeper understanding. Recognizing VR's limitations highlights the need for a comprehensive approach [21]. Tools like interactive 3D models and dynamic visualizations enable users to explore designs, identify issues, and experiment with solutions. This holistic approach is essential for meeting the multifaceted comprehension requirements of highly complex injection mold designs.

Several researchers have made contributions to understanding and applying virtual reality in various industrial contexts. Their studies explore the use of virtual prototypes for various applications, such as product communication, feature evaluation and factory planning. For instance, Bordegoni et al. [22] investigated the use of virtual prototypes, incorporating visual, haptic, and acoustic elements, as an alternative to physical prototypes for communicating new products and assessing their features. Ferrise et al. [23] conducted case studies on interactive virtual prototypes that replace physical models in conceptualization and product design activities. In the context of factory planning, Gebhard et al. [24] developed a VR application that facilitates factory design through realistic tours in digital models. On the other hand, Sampaio et al. [25] focused on creating virtual models to support decision-making in construction management and maintenance. Lastly, Boton [26] proposed the inclusion of the construction process in 4D models for constructability analysis. By contextualizing these studies within the broader landscape of VR applications, it becomes evident that VR technology holds significant promise for various industrial processes, including those within the injection molding sector. However, despite these advancements, there remains a notable gap in research regarding the direct application of VR to the geometric and functional validation of complex injection molds. To date, the only recognized study in the domain of injection molding is the research conducted by Sun et al. [27]. Their work concentrates on developing a virtual training simulation platform for injection molding machines, aimed at enhancing safety through comprehensive training programs. Unfortunately, this study focuses only on training for injection machines, neglecting considerations related to mold tooling.

Given the highlighted gaps, there is a clear absence of research focusing on the comprehensive validation of design injection molds, surpassing the limitations of conventional CAD systems. The influence of specific parameters on the validation processes, assessed through quantitative statistical analysis with expert involvement, remains largely unexplored. This presents an area with significant potential for advancing understanding and refining injection mold manufacturing processes, especially given the anticipated impact of this technology in the future. This underscores the need for further exploration and innovation, particularly in leveraging VR to address the complex challenges inherent in mold design validation and production efficiency. Despite the existing literature and reviews, significant research gaps persist regarding the application of VR in industrial contexts [28,29]. The research presented in this paper surpasses the state of the art aiming to bridge this gap by evaluating the effectiveness of multidimensional CAD-virtual validation methods compared to traditional CAD systems. Through animated VR representations, the research enhances comprehension and validation of mold functionality, demonstrating its capability to streamline processes effectively. The objective is to showcase the efficiency enhancements achieved through VR devices in mold design validation, compared to conventional CAD methods. Ultimately, the goal is to illustrate how VR technology can substantially reduce validation times and manufacturing errors in injection molding, thereby offering transformative advantages for manufacturers and emerging as an essential tool in the industry.

## 2. Materials and Methods

This research work focuses on addressing the application of novel VR techniques and tools, compared to the traditional use of CAD tools in the field and industrial sector of plastic injection molds. Through this research, we delineate the disparities and parallels, along with the advantages and drawbacks, arising from the comparison between traditional techniques used in the industrial sector of injection molds and the application of new VR tools. The presented work specifically focuses on significant areas such as CAD geometric design, comprehensive analysis of assembly functionality, detection of design flaws and component interferences, and validation of assembly and disassembly operations within injection mold processes. This section provides a complete description of the implementation and development process of the proposed methodology.

Presently, there exists an extensive array of commercial software designed to integrate and operate with VR devices and technology. However, their compatibility and integration with CAD-type tools remain limited. Generally, VR tools lack support for CAD file formats and offer restricted interaction with 2D and 3D geometric entities and operations. That is why, for the methodology proposed in this manuscript, the free software Blender 4.0 is used as an intermediate resource to the commercial software Aurora VR [30] (see Figure 1) as a link between the VR environment and operations and the geometric entities modelled in the CAD environment. Aurora VR does not allow importing industrial components and assemblies in CAD format into the VR environment directly. With this in mind, a series of operations are available that allow the user to interact with the imported geometry, including the following: cutting plane, rotation/translation movements in space, scaling, annotation, measurement and animation (refer to Figure 2). On the other hand, for the 3D CAD modelling of the components analyzed in this manuscript, the commercial software CATIA V5 R21 [31] has been used. The VR integration in this methodology employs HP Reverb G2 VR glasses and a VR-ready GPU-equipped PC (Core i5, 3.60 GHz, and 16 GB RAM), as illustrated in Figure 3.

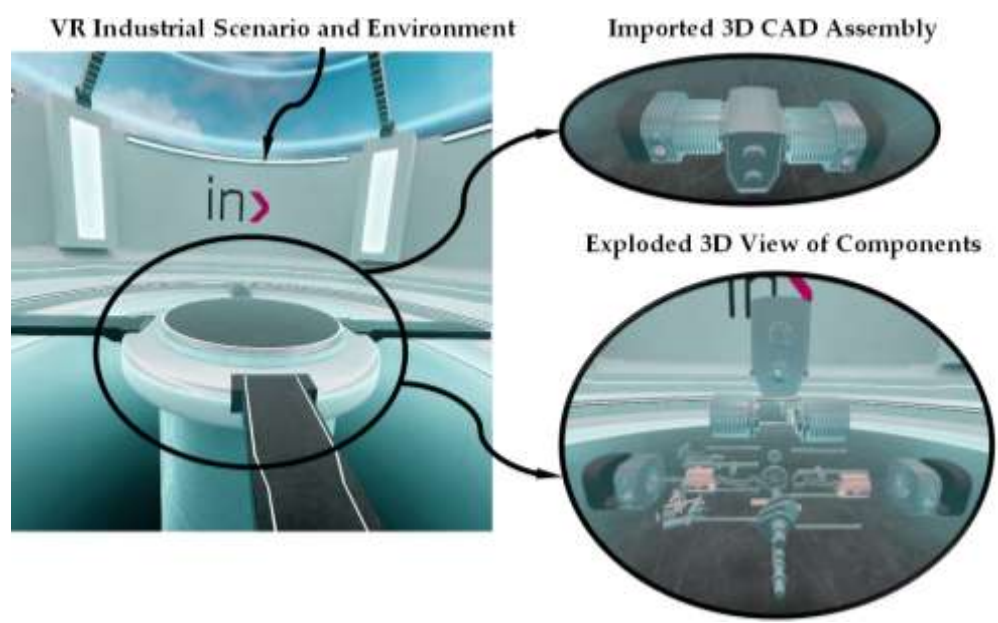

**Figure 1.** Definition of the VR software used.

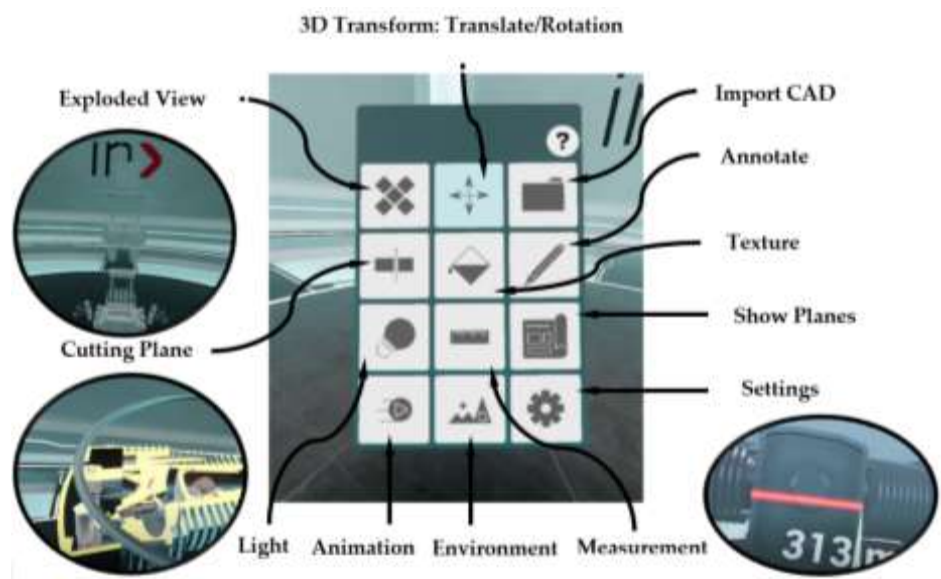

**Figure 2.** VR Aurora work environment.

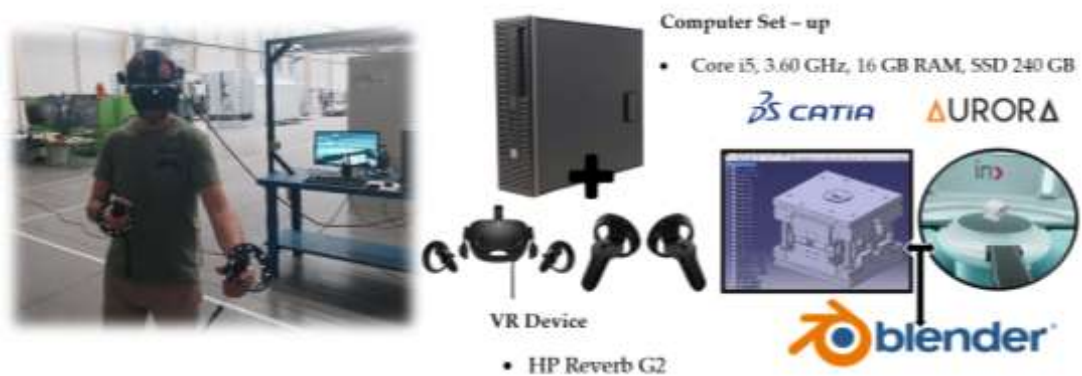

**Figure 3.** CAD—VR devices and environment.

## 2.1. Three-Dimensional CAD Modelling

The advancement of commercial 3D CAD modeling software has progressed alongside the evolving demands of various engineering sectors. However, the integration of VR technology, which is currently experiencing rapid growth, with 3D CAD geometric modeling processes still faces significant technical challenges. One of the primary obstacles is the format incompatibility between CAD geometric files and their integration into VR environments or software. To avoid this, all the 3D CAD models analyzed in this research were transformed from CATIA V5 R21 to STL format (STereoLithography—Standard Triangle Language). This STL meshing or faceted format allows the topology of the components to be preserved without influencing the precision or geometric details of the components. Next, in Blender 4.0, we will export to \*.GLB format adapted and compatible with the Aurora VR Software tools (Figure 3).

The methodology detailed in this manuscript aims to establish a multidimensional environment utilizing interactive VR technology. Its objective is to enhance and streamline the recognition and comprehension of injection mold design through various visualization tools. These tools assist in detecting interferences between components, identifying design errors, and providing support for assembly, disassembly, and maintenance operations.

To accomplish this objective, this manuscript presents 3D CAD models of three injection molds used in industrial plastic part manufacturing for the automotive industry. As shown in Figure 4, these molds display design errors and interferences among their components. In the industrial sector of injection molds, maintaining a safe distance between all mold components and systems is a critical design requirement. For example, complex elements of the cooling system often fail to meet the minimum safety distance, or elements of the ejection system may interfere with certain main mold plates.

Conducting these inspections in a CAD environment is laborious and time-consuming for designers. Moreover, the complexity of conventional injection mold mechanisms and components is substantial, and 3D CAD modeling environments lack real-time visualization and interaction operations with three-dimensional models, hindering comprehensive design understanding. Figure 4 shows examples of common errors made during the design of the main elements and systems of the injection mold. Two main design errors were identified. On the one hand, the minimum safety distance, Dsafety (see Figure 4), is not maintained between the components of the injection mold. Typically, the industrial sector of plastic injection molds establishes a safety distance, Dsafety (see Figure 4), which is equal to 10 mm, to maintain the structural integrity of the injection mold. On the other hand, interference may occur between them. Hence, this research evaluates whether VR technology can significantly advance the industrial sector of injection molds.

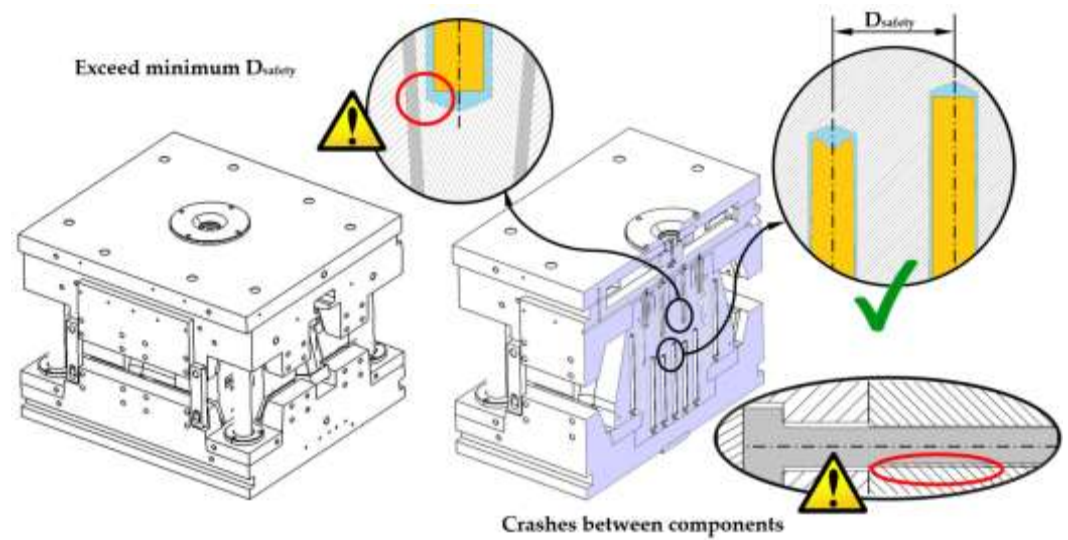

**Figure 4.** Injection mold 3D CAD model. Design error detection and interferences.

## 2.2. VR Interaction and Integration

Traditionally, designers interact with 3D CAD models within the CAD software environment using its built-in visualization operations, along with the available PC hardware. Despite efforts by CAD application developers to enhance this interaction using traditional tools, technical obstacles persist, impeding real-time manipulation and visualization of designs. These limitations directly impact geometry review operations. Particularly when assessing complex geometric assemblies with injection molds, the interaction becomes challenging. To enhance and streamline this process, this methodology proposes leveraging the advantages and versatility offered by VR technology and devices for real-time analysis and inspection operations of injection molds. Additionally, it aims to evaluate the extent of improvement and satisfaction that this technology can bring to the industrial sector.

Firstly, as illustrated in Figures 1–3, the software utilized as the VR interaction environment is Aurora VR (Invelon, Ingroup, Gipuzkoa, Spain) [30], paired with HP Reverb G2 VR glasses. The HP Reverb glasses feature two 2.89-inch LCDs with a resolution of 2160 × 2160, DIP adjustment, a 90 Hz frequency, and four integrated tracking cameras. The VR device calibration is configured at 60 mm and operates at a frequency of 90 Hz.

The integration of the VR device is straightforward and does not require an additional configuration process. The Aurora VR software automatically detects the device and sets it up for immediate use. Regarding the import of the geometry of the analyzed injection molds, it is performed using the STEP-neutral CAD compatibility format.

The VR environment prioritizes essential factors such as immersion, visual quality, and smooth real-time interaction. The primary objective of the workspace in the VR environment is to attain adaptability, enabling the understanding of geometrically complex sets like injection molds. This immersive environment allows for a detailed evaluation of the scale, proportions, and layout of the design within three-dimensional space.

As shown in Figure 2, once the CAD model is imported, the system detects the user's movements and facilitates interaction with the VR environment and the imported models through the operations palette. Utilizing the operations palette does not require programming skills. Users can execute various operations by simply selecting icons, including the following:

- Separate: This function enables users to adjust the positioning of components within three-dimensional space. Typically, it is utilized to interact with assemblies comprising diverse components that share positional relationships. As demonstrated in Figure 5, it facilitates the step-by-step disassembly of the assembly. This simulated and virtual disassembly process enables access to interior and inaccessible parts of the assembly, offering significant utility and potential application in plastic injection molds. Indeed, these components are typically heavy and impractical to manipulate individually without auxiliary mechanisms. Moreover, this operation aids in planning the assembly and disassembly of the entire injection mold, adapting it to the real workspace and the injection machine. Furthermore, it enhances the understanding of the layout and position of each mold component, facilitating maintenance or repair tasks. Lastly, after any modification, selecting this operation again restores the assembly components to their original positions.

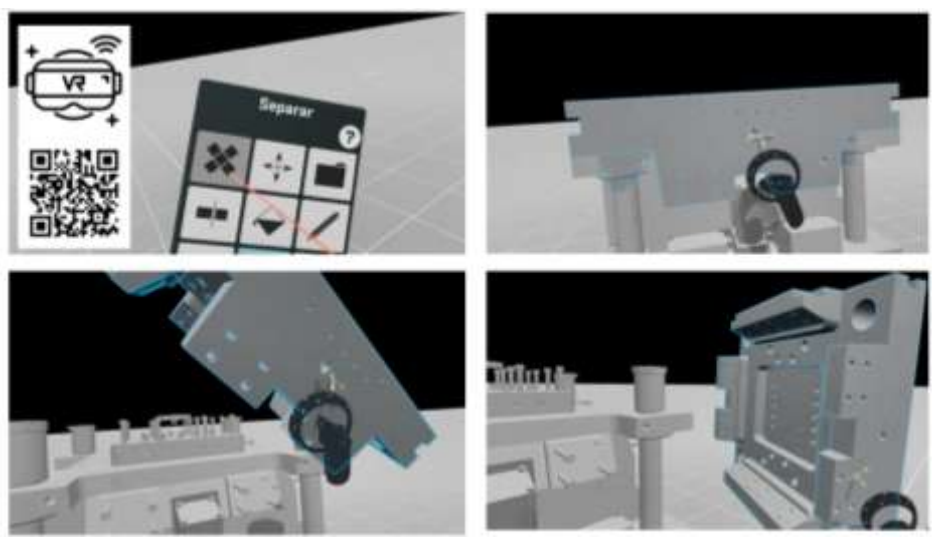

**Figure 5.** Separation and disassembly operation of assemblies applied to injection molds.

- Cutting plane: this operation, as shown in Figure 6, allows cutting from a plane defined, interactively, by the user. In this way, each component of the assembly can be accessed and viewed, maintaining its assembly position. In particular, this operation is useful to perform an analysis of interferences between components or detect if there is any design or assembly error between components.

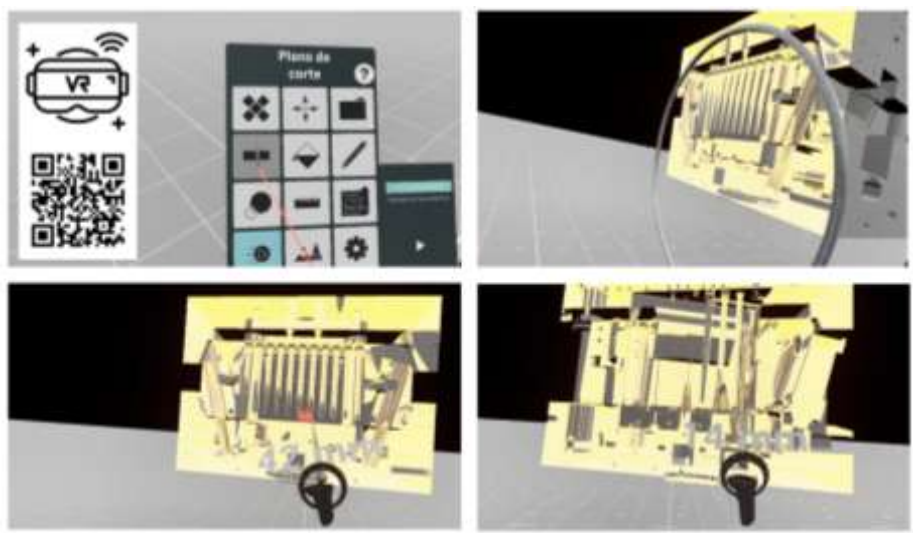

**Figure 6.** Plane cutting operation applied to injection molds.

- Three-dimensional transform: through this operation, the geometry can be oriented through translation movements along the three main axes of three-dimensional space, rotation movements about these axes and movements with all degrees of freedom. Likewise, three-dimensional models can also be uniformly scaled interactively.
- Annotate: This operation, as shown in Figure 7, allows the annotation and addition of information in the VR environment.

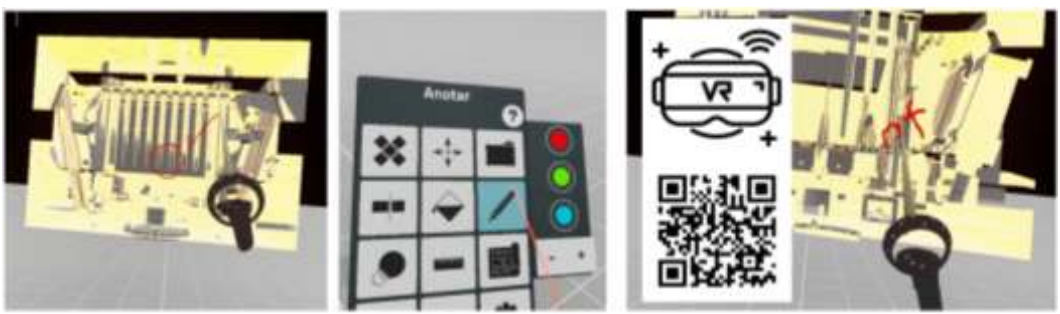

**Figure 7.** Plane cutting operation applied to injection molds.

- Measurements: As shown in Figures 6 and 8, this tool empowers users to directly measure the spatial separation between selected points, defined via controllers associated with the VR device. Consequently, the accuracy of the obtained measurements is hinged on the precision with which users specify the initial and final points on the analysis geometry. When executed correctly, without deviations from the proposed geometric elements, the resulting measurements are devoid of errors. Moreover, it is worth underscoring the virtual measurement's sensitivity, which operates within the millimeter magnitude range. This sensitivity aligns with typical design errors encountered during the modeling process of primary injection mold systems.

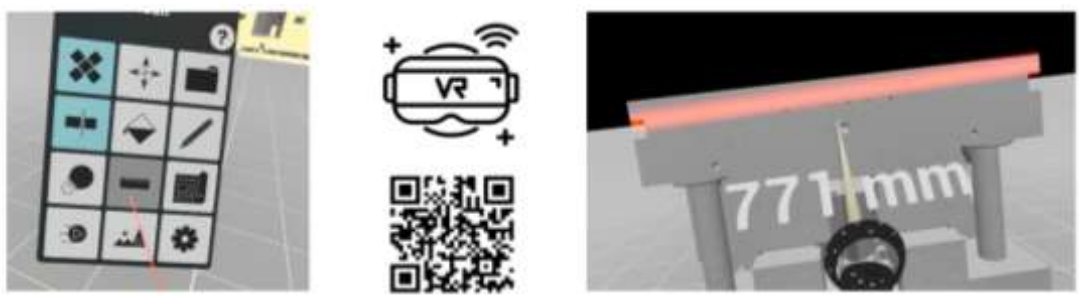

**Figure 8.** Measurement operation applied to injection molds.

- Animation: View animations or motion studies are defined before importing the model into the VR environment. From this operation, the simulation of the breakdown of an assembly or the representation of an exploded view of all its components can be programmed (see Figure 9).

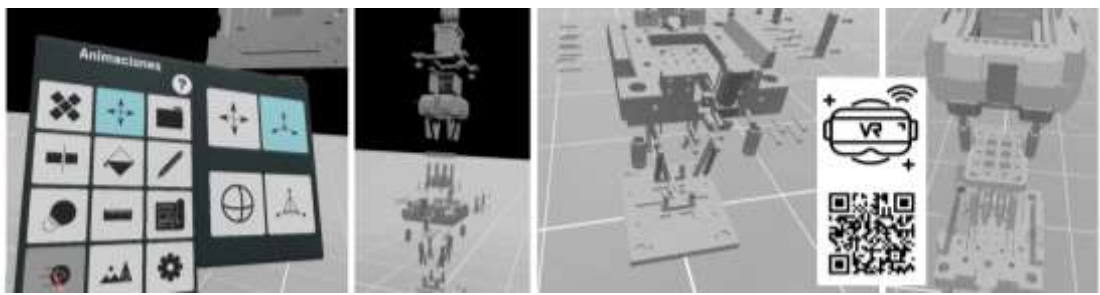

**Figure 9.** Exploded view of the geometry of an injection mold in VR.

- Displacements: Through this operation, the user can move to any position in the VR environment (see Figure 10).

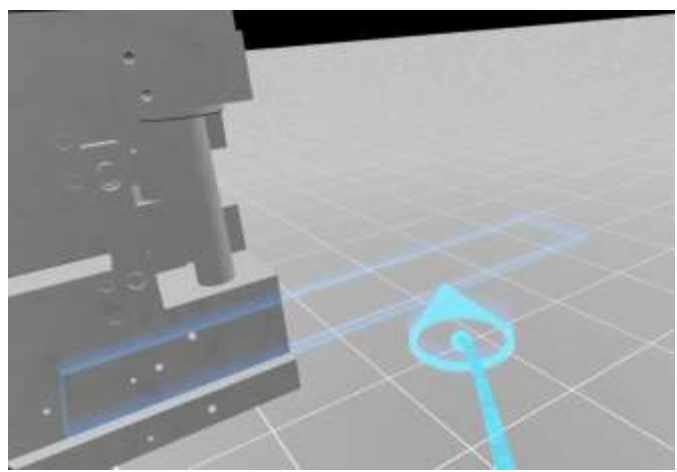

**Figure 10.** User movement operation throughout the VR environment.

## 2.3. Application of VR Technology to an Industrial Activity

After outlining the VR environment and user interaction capabilities concerning injection molds, the following section details the proposed activity aimed at evaluating this technology within a realistic industrial context.

To achieve this goal, two injection mold design evaluation sessions were conducted in collaboration with an engineering company highly specialized in the industrial and automotive sectors (refer to Figure 11). Renowned for its expertise in designing highly complex plastic parts, this company excels in crafting advanced injection molds, produced both in-house and through external suppliers via subcontracting.

Within the company, the mold design review sessions were conducted by injection mold manufacturing experts with vast experience in the sector. These meetings serve as validations of the work conducted by the mold design development team and provide recommendations for enhancing the design. Furthermore, they aim to impart specialized information to suppliers, facilitating efficient decision-making whenever feasible.

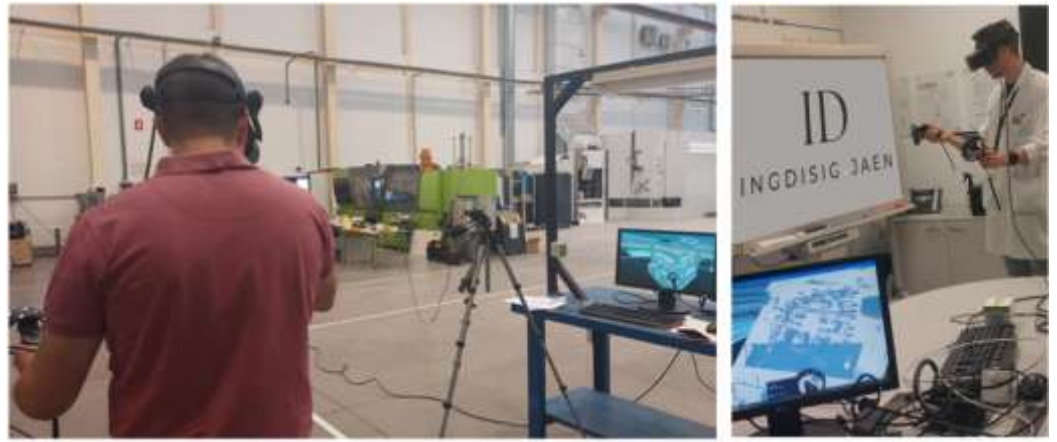

**Figure 11.** Development of injection mold design verification activities with VR technology.

The sessions were conducted in small groups, with no more than three participants, who reviewed the designs and shared their experiences from various perspectives.

In the initial phase, an experienced mold maker specializing in CAD design within the company selected three complex injection molds for evaluation. The injection molds (refer to Figure 12) were chosen based on design, size, and complexity criteria.

Subsequently, a list of eight common defects in the design of injection molds was compiled. These defects were derived from real situations and the extensive experience of the company. It is noteworthy that although the defects may be similar in different case studies, they rarely coincide in the same location of the mold. The errors or defects in the presented molds were classified into the following groups:

- Sizing errors in the cooling and ejection system;
- Accessibility problems for assembling the parts in the mold;
- Interferences between the different parts of the mold;
- Deficiencies in the logic of the cooling layout;
- Absence of parts or poor placement of small elements, such as sealing gaskets;
- Misalignment of mold components;
- Lack of robustness in critical areas of the mold.

These defects represent key areas that require attention and correction during the injection mold design and manufacturing process.

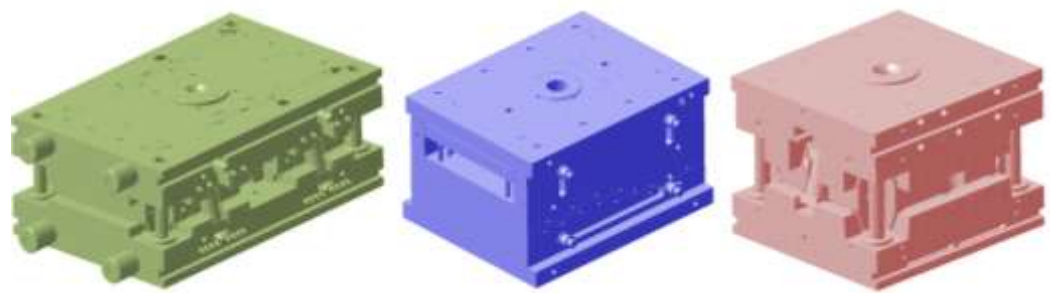

**Figure 12.** Injection molds used for the development of the activity.

A total of 22 participants were involved in the development of this methodology, distributed across 7 work teams. These review teams were formed in a heterogeneous manner, ensuring that each member contributed a different vision, based on her professional experience. For each group, a period of 40 min was assigned to complete the activity in its entirety. During this, the two design review approaches were explored on various 3D CAD models. The process of the proposed activity and its evaluation was developed as follows:

- Introduction: 15 min.
- Design review using VR or Catia V5 support in 3D model.
- Second review: 15 min.
- Design review using the alternative method (not used in the previous step) on the 3D model.
- Quiz: 10 min (see Table 1).

In the introduction phase, the procedural details were outlined, and participants were prompted to identify potential design flaws in the 3D molds provided. Teams initiated the review process by utilizing the VR tools as described earlier (refer to Figure 3). All participants utilized the HP Reverb G2 VR devices. While one team member engaged with the VR devices for the review, the remaining team members observed the transmitted images from the VR software on a designated monitor (refer to Figure 11). Following this, after 20 min, the group shifted to utilizing CAD modeling software to examine the second injection mold model. During this phase, one team member operated the CAD modeling software while the remaining members observed the screen. Similar to the prior examination, all individuals within each group engaged with the CAD modeling software. To offset any potential learning bias, the activity alternated the start between the VR and CAD phases among successive groups. Upon the identification of a defect by a group member and subsequent agreement by the team, a note documenting the identified issue was recorded. Following this, each participant filled out a questionnaire assessing both mold design review sessions using a Likert scale ranging from 1 to 7 points (refer to Table 1).

**Table 1.** Survey carried out to assess the application of VR technology in the injection mold industrial sector.

| Questions | Likert Scale—Limits Value | |
|---|---|---|
| Detection of dimensional errors | Q1. Evaluate the conventional CAD software environment in operations to identify dimensional design errors in complex injection molds | 1—Very dissatisfied 7—Very satisfied |
| | Q2. Evaluate the VR environment in operations to identify dimensional design errors in complex injection molds | 1—Very dissatisfied 7—Very satisfied |
| Functional error detection | Q3. Evaluate the level of clarity and ease of acquiring information provided by the conventional CAD software environment during the functional error review of injection mold design. | 1—Very dissatisfied 7—Very satisfied |
| | Q4. Evaluate the level of clarity and ease of acquiring information provided by the VR environment during the functional error review of injection mold design. | 1—Very dissatisfied 7—Very satisfied |
| Design understanding | Q5. Evaluate the level of efficiency of the conventional CAD software environment to help understand the concept or detailed design of a complex injection mold | 1—Very dissatisfied 7—Very satisfied |
| | Q6. Evaluate the level of efficiency of the VR environment to help understand the concept or detailed design of a complex injection mold | 1—Very dissatisfied 7—Very satisfied |
| Assembly and disassembly validation | Q7. Evaluate the help offered by the conventional CAD software environment in validation operations for the assembly and disassembly of the injection mold design. | 1—Very dissatisfied 7—Very satisfied |
| | Q8. Evaluate the help offered by the VR environment in validation operations for the assembly and disassembly of the injection mold design | 1—Very dissatisfied 7—Very satisfied |

| Influence on the training process | Q9. Evaluate the influence of the conventional CAD software environment on the training process of new professionals, without prior experience, in this industrial sector. | 1—Very dissatisfied<br>7—Very satisfied |
|---|---|---|
| Influence on the training process | Q10. Evaluate the influence of the VR environment in the training process of new professionals, without previous experience, in this industrial sector. | 1—Very dissatisfied<br>7—Very satisfied |
| Satisfaction degree | Q11. Evaluate the level of satisfaction with the use of VR technology in the field of injection mold design | 1—Very dissatisfied<br>7—Very satisfied |
| | Q12. Evaluate the level of satisfaction with the use of VR technology in the assembly and disassembly of injection molds | 1—Very dissatisfied<br>7—Very satisfied |

### 2.4. Participants

As previously outlined, the activity conducted in this research involves expert professionals with extensive experience in the design and manufacture of injection molds. Specifically, as depicted in Table 2, these professionals have been categorized into three groups: engineers, mold technicians, and designers. Additionally, relevant aspects and prior knowledge pertinent to this research have been defined for each participant. These aspects include the level of familiarity with the use of VR technology, the experience level in CAD design within the injection mold field, and expertise in assembly and disassembly. All of these aspects are assessed using a Likert scale, with extreme values ranging from 1 to 7. Table 2 provides a summary of the main aspects characterizing the participants in the proposed activity.

**Table 2.** Sample demographics.

| Parameters | Number of People | % |
|---|---|---|
| Professional profile | | |
| Engineer | 14 | 63.6 |
| Mold technician | 4 | 18.2 |
| Designer | 4 | 18.2 |
| Previous VR technology experience | | |
| 1–2 | 19 | 86.4 |
| 3–4 | 3 | 13.6 |
| 5–6 | 0 | 0.0 |
| 7 | 0 | 0.0 |
| Previous CAD design experience | | |
| 1–2 | 2 | 9.1 |
| 3–4 | 8 | 36.4 |
| 5–6 | 7 | 31.8 |
| 7 | 5 | 22.7 |
| Previous experience in mold assembly and disassembly | | |
| 1–2 | 5 | 22.7 |
| 3–4 | 5 | 22.7 |
| 5–6 | 9 | 40.9 |
| 7 | 3 | 13.6 |

As evident, the majority of participants lack previous experience in the use of VR technology, with 84.6% considering themselves novices in this regard. Conversely, over 50% possess significant prior knowledge in designing complex injection molds using CAD-type software, along with previous experience in their assembly and disassembly processes. Therefore, the selected sample population for conducting the proposed activity aligns with the requirements and premises established to evaluate the influence and impact of virtual reality technology on the plastic industrial injection molding sector.

To ensure the adequacy of the selected population size for subsequent statistical analysis of the conducted survey, an initial population calculation is performed using GPower statistical calculation software (version 3.1, University of Düsseldorf, Düsseldorf, Germany) [32,33]. Considering that the survey's objective, included in the conducted activity, involves comparing and evaluating the influence of VR technology versus conventional CAD environments, a *t*-test comparative analysis is proposed. Table 3 outlines the magnitude and selection of parameters established in the GPower statistical calculation software for calculating the initial population.

**Table 3.** Parameters defined for the statistical calculation of the initial population.

| Input Parameters | Selection |
|---|---|
| Test family | *t*-test |
| Statistical test | Means: Difference between two dependent means (matched pairs) |
| Type of power analysis | A priori: Compute required simple size |
| Tail | Two |
| Effect size dz | 0.80 |
| $\alpha$ err prob | 0.05 |
| Power (1−$\beta$ err prob) | 0.95 |

A statistical test was used to establish the type of statistical test to perform. In particular, for the present statistical analysis, the "Difference between two dependent means" was used since it is a very common option to compare the average between two independent groups. The type of power analysis represents the statistical power of the research. Since the calculation of the initial population has been carried out prior to carrying out the experiment, the a priori option has been defined. Tail determines whether the type of contrast is unilateral or bilateral, effect size dz determines the size of the effect in the hypothesis made that may exist in our research work, and $\alpha$ err prob/power (1−$\beta$err prob) represents the power statistics or confidence level for our research work. In this way, a bilateral contrast has been defined given that, a priori, the hypothesis is raised that the results obtained in the evaluation of VR technology will be different from the results in the evaluation of conventional CAD environments. A value of 0.65 has been defined as the effect size, given that the effect size is considered between medium (0.5) and large (0.8). Furthermore, for the statistical power or confidence level, a standard statistical magnitude equal to 0.05 and 0.95 has been determined.

Finally, Table 4 shows the results obtained from the statistical analysis of the initial population carried out. As can be seen, the minimum size of the initial population proposed for activity is 21 individuals. This, in turn, is lower than the size of the population finally established, equal to 22. Therefore, in this way, the size of the population established for the activity proposed in this research work is validated. Moreover, other parameters are essential to ensure the validity and reliability of the research findings. The noncentrality parameter ($\delta$) indicates the effect size, representing the difference between the null hypothesis and the true population mean. In our analysis, $\delta$ was calculated as 3.8951, indicating a substantial effect size.

The critical *t*-value (2.0860) was used as the threshold for determining statistical significance. This helps assess whether the observed differences are likely due to chance or represent true differences in the population. With a degree of freedom (Df) of 20, which accounts for the variability in the sample, and a total sample size of 21, our study had adequate statistical power to detect meaningful effects. The calculated actual power of 0.9592 indicates the likelihood of detecting the true effects when they exist.

**Table 4.** Parameters defined for the statistical calculation of the initial population.

| Output Parameters | Selection |
|---|---|
| Non-centrality parameter δ | 3.8951 |
| Critical t | 2.0860 |
| Df | 20 |
| Total simple size | 21 |
| Actual power | 0.9592 |

## 3. Results and Discussion

The statistical analysis of the results obtained from the proposed survey (see Table 1) is carried out using the statistical calculation program IBM SPSS Statistics (version 27, IBM, Armonk, NY, USA) [34]. Firstly, to evaluate the validity of the questions proposed for the survey, Cronbach's Alpha test was used. Cronbach's Alpha test evaluates the interrelationship and reliability of the responses, indicating the internal consistency and coherence of the questions and Likert scales defined in the survey carried out. As shown in Table 5, the value obtained using Cronbach's Alpha test equals 0.749. This result confirms that the defined Likert scale and the global structure of the survey, proposed in the activity carried out, are consistent. Cronbach's Alpha result is within the limit values for this test, that is, between 0.70 and 0.9. This validation is important to reinforce the reliability of the statistical analysis and guarantee that the results serve as a solid basis for subsequent interpretation. Table 5 shows the results obtained from Cronbach's Alpha test.

**Table 5.** Cronbach's Alpha test results.

| Cronbach's Alpha | Cronbach's Alpha Based on Standardized Items | N° of Likert Questions |
|---|---|---|
| 0.749 | 0.749 | 12 |

Next, before proceeding to carry out any operation or statistical comparison test, it is necessary to determine whether the distribution of the responses obtained, for the survey carried out, maintains a normal distribution. This normality analysis determines the type of comparison tests or *t*-tests that best fit the distribution presented by the data obtained. The normality test used in this research work is the Shapiro–Wilk test, given that the population size is less than 50. Table 6 shows the results obtained from the normality test.

**Table 6.** Normality test results for the proposed survey questions.

| Questions | | *p*-Values |
|---|---|---|
| Detection of dimensional errors | Q1—CAD | 0.042 |
| | Q2—RV | 0.000 |
| Functional error detection | Q3—CAD | 0.002 |
| | Q4—RV | 0.000 |
| Design understanding | Q5—CAD | 0.045 |
| | Q6—RV | 0.000 |
| Assembly and disassembly validation | Q7—CAD | 0.048 |
| | Q8—RV | 0.000 |
| Influence on the training process | Q1—CAD | 0.078 |
| | Q2—RV | 0.000 |
| Satisfaction degree | Q11 | 0.002 |
| | Q12 | 0.000 |

As can be seen, the results obtained from the normality test (see Table 6) show that the data obtained from the survey carried out do not fit a normal distribution, with a level of significance or confidence *p*-value < 0.05. Therefore, it is necessary to use non-parametric tests before carrying out the subsequent statistical tests compared to the *t*-test.

Next, based on the results obtained, contextualized in the proposed Likert scale, a descriptive statistical analysis is carried out for each of the questions defined in the survey carried out (see Table 1). This analysis provides, based on the mean, standard deviation and variance values and information about the central tendency and dispersion of the data obtained. Table 7 shows the results obtained for this descriptive statistical analysis.

**Table 7.** Results obtained for the descriptive statistical analysis.

| Questions | | Mean | Standard. Deviation | Variance |
|---|---|---|---|---|
| Detection of dimensional errors | Q1—CAD | 4.86 | 1.207 | 1.457 |
| | Q2—RV | 6.09 | 0.610 | 0.372 |
| Functional error detection | Q3—CAD | 5.36 | 1.217 | 1.481 |
| | Q4—RV | 6.27 | 0.883 | 0.780 |
| Design understanding | Q5—CAD | 5.14 | 0.990 | 0.980 |
| | Q6—RV | 6.36 | 0.727 | 0.529 |
| Assembly and disassembly validation | Q7—CAD | 5.00 | 0.926 | 0.857 |
| | Q8—RV | 6.73 | 0.550 | 0.303 |
| Influence on the training process | Q1—CAD | 4.82 | 1.053 | 1.109 |
| | Q2—RV | 6.64 | 0.658 | 0.433 |
| Satisfaction degree | Q11 | 6.09 | 0.921 | 0.848 |
| | Q12 | 6.59 | 0.666 | 0.444 |

Likewise, as shown in Figures 13–15, an analysis of relative frequencies has been carried out for each of the questions evaluated using a Likert scale for the proposed survey. Figure 13 shows the comparison of results obtained when detecting dimensional design errors using VR technology and conventional CAD environments (see Table 1). As can be seen, the results of the survey show that dimensional errors, such as poor sizing of injection mold systems or interference detection, are detected more efficiently using VR technology. A significant majority, 86.37%, positively evaluated the application of this technology in detecting dimensional errors, scoring above 6 on the Likert scale. Likewise, Figure 13 also compares the results obtained when detecting design functional errors using VR technology and conventional CAD environments (see Table 1). Once more, the survey results indicate that functional errors, such as inadequate layouts of injection mold systems or the absence of fundamental components, are more efficiently detected using VR technology. In particular, 59.09% have positively valued, scoring above 6 on the Likert scale, the application of this technology in the detection of functional errors. Similarly, the trend in these assessments is maintained for the rest of the questions asked (see Figures 14 and 15). However, in areas related to the understanding of the design and functionality of injection molds, the validation of their assembly and disassembly tasks and the influence on the learning and training process for future professions in this sector, the evaluations obtained by applying VR technology stand out significantly compared to the use of conventional CAD modelling environments. In particular, 50%, 77.27% and 72.73% of the mold experts value the application of VR technology for the understanding of the design and functionality of injection molds, the validation of assembly and disassembly tasks and the influence on the learning and training process for future professionals in this sector, respectively, very positively, with a score of 7 on the Likert scale.

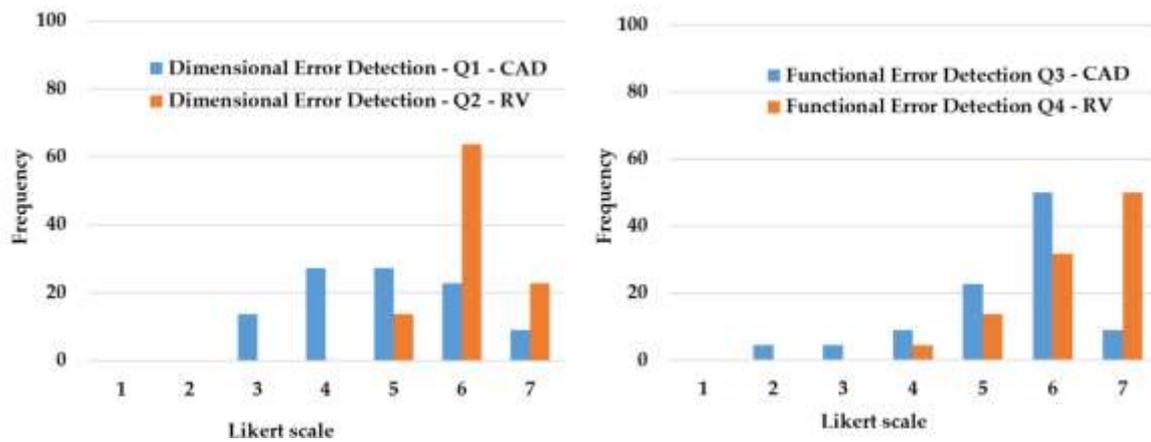

**Figure 13.** Response frequency obtained for the error detection questions.

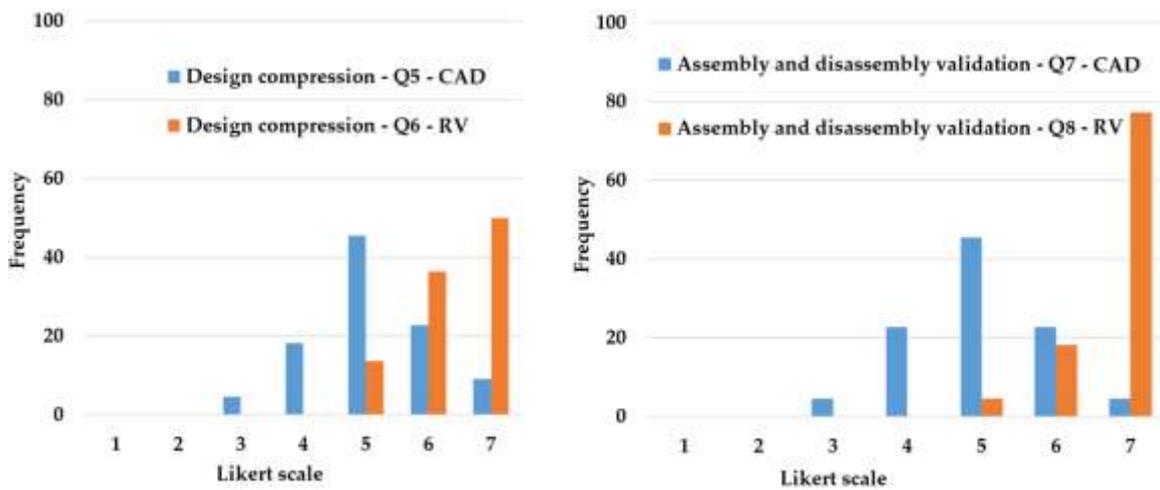

**Figure 14.** Response frequency obtained for design compression and assembly and disassembly validation questions.

Finally, as depicted in Figure 15, the satisfaction level regarding the application of VR technology in the injection molding industrial sector is assessed through two questions. Firstly, the satisfaction level with the utilization of VR technology in injection mold design processes is evaluated. It was observed that molding experts, with a cumulative relative frequency of 72.72%, positively appraise the application of VR technology in injection mold design processes, indicating a Likert scale score. Furthermore, the satisfaction level with the application of VR technology in assembly and disassembly tasks of injection molds is also assessed. It is evident that mold experts, with a relative frequency of 90.91%, highly appreciate the application of VR technology in these tasks, indicating a Likert scale score greater than 6.

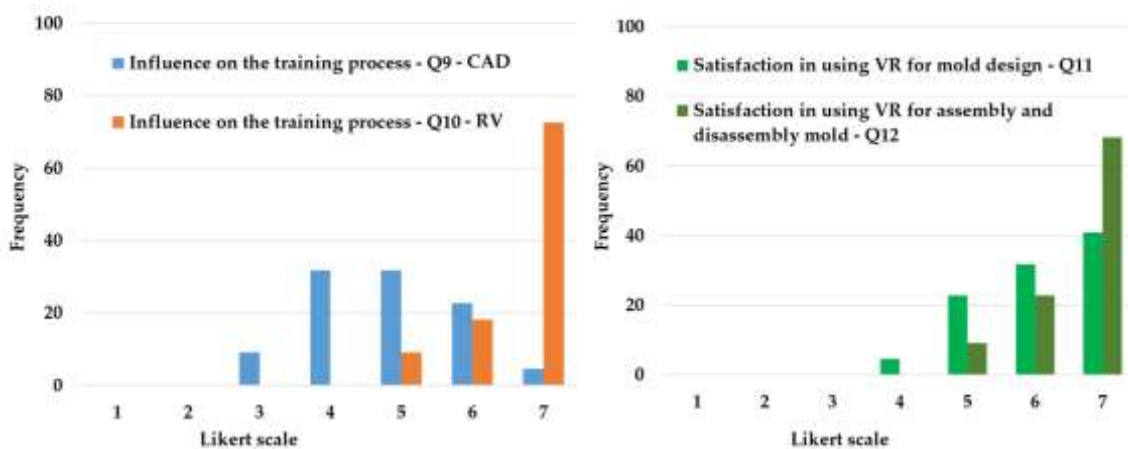

**Figure 15.** Response frequency obtained for the influence on training process and the satisfaction in the use of VR technology questions.

Next, to analyze the impact associated with the application of VR technology in the injection molding industrial sector, it is proposed to perform a non-parametric comparison test, the W Wilcoxon signed-rank test. This non-parametric test examines whether the utilization of VR technology has significantly influenced the outcomes concerning the detection of dimensional and functional design errors, understanding of design and functionality, and validation of assembly tasks and disassembly, as well as its impact on the learning and training process for future professionals. Table 8 shows the results of the W Wilcoxon signed-rank test.

**Table 8.** W Wilcoxon signed-rank test results.

| Questions | | *p*-Values |
|---|---|---|
| Detection of dimensional errors | Q1—CAD<br>Q2—RV | 0.002 |
| Functional error detection | Q3—CAD<br>Q4—RV | 0.010 |
| Design understanding | Q5—CAD<br>Q6—RV | 0.002 |
| Assembly and disassembly validation | Q7—CAD<br>Q8—RV | 0.001 |
| Influence on the training process | Q9—CAD<br>Q10—RV | 0.001 |

The results of the W Wilcoxon signed-rank test, presented in Table 8, show that there are significant differences, with a degree of significance or confidence level *p*–value < 0.05, between the results obtained for the activity tasks carried out through conventional CAD environments and the application of VR technology. In this way, through this comparative test, it is confirmed that the application of VR technology has a direct and significant impact on the industrial sector of injection molds, especially in the areas of detection of dimensional and functional design errors, understanding of the design and functionality validation of assembly and disassembly tasks and influence on the learning and training process for future professionals.

Finally, to conclude the statistical analysis of the results obtained, a non-parametric test, the Kruskal–Wallis test, was carried out to check whether the professional figure of the mold experts who have participated in the proposed activity and their level of familiarity with the use before VR technology significantly influences the evaluation of this as a useful application tool in the injection molding industrial sector. Table 9 shows the results obtained for the non-parametric Kruskal–Wallis test. As can be seen, there are no significant differences in the results obtained with a level of confidence or significance *p*-value > 0.05. That is, the professional figure of the molding experts and the level of familiarity with the previous use of VR technology do not significantly influence the assessment of the application of VR technology in the industrial sector of injection molds.

**Table 9.** Kruskal–Wallis test results.

| Questions | | Influence on Professional Figure *p*-Values | Influence on Previous Use of VR Technology *p*-Values |
|---|---|---|---|
| Detection of dimensional errors | Q1—CAD Q2—RV | 0.365 | 0.668 |
| Functional error detection | Q3—CAD Q4—RV | 0.964 | 0.882 |
| Design understanding | Q5—CAD Q6—RV | 0.398 | 0.291 |
| Assembly and disassembly validation | Q7—CAD Q8—RV | 0.473 | 0.374 |
| Influence on the training process | Q9—CAD Q10—RV | 0.744 | 0.719 |

## 4. Conclusions

This paper introduces an innovative multimodal virtual environment designed specifically for validating mold design complexities, integrating interactive 3D models and dynamic visualizations to facilitate exploration of detailed mold designs. This research represents a significant advancement beyond current practices, aiming to address existing gaps by evaluating multidimensional CAD-virtual validation methods against traditional CAD systems. Through the utilization of animated VR representations, the study significantly enhances the comprehension and validation of mold functionality, ultimately streamlining various processes involved. Furthermore, the research underscores the disparities, parallels, and advantages and drawbacks between traditional techniques and VR applications in injection mold design. This comprehensive analysis provides valuable insights into the potential of VR technology to revolutionize mold validation processes and reshape the landscape of industrial design practices.

The survey findings reveal that dimensional errors are effectively identified using VR technology, with 86.37% of respondents providing positive evaluations of its application. Similarly, functional errors, such as inadequate layouts, are efficiently detected with VR, garnering positive feedback from 59.09% of participants. Mold experts place high value on VR technology for comprehending design functionality (50%), validating assembly and disassembly tasks (77.27%), and facilitating the learning process (72.73%). Statistical analyses, including the Wilcoxon test, uncover significant differences in interference detection, topology tracking, and validation of assembly and disassembly accessibility between traditional CAD and VR methods. These results underscore VR's potential to enhance mold design validation processes, with a relative frequency analysis demonstrating efficiency gains in mold tracking and validation using VR devices. Furthermore, the non-parametric Kruskal–Wallis test indicates that professionals' familiarity with VR technology does not significantly influence its assessment in the injection molding industry.

This research highlights the crucial role of VR technology in revolutionizing injection molding design processes. Its ability to enhance efficiency and accuracy signifies a significant advancement, facilitating streamlined operations and improved outcomes in the manufacturing sector. This transformative potential offers manufacturers considerable advantages, firmly establishing VR as an indispensable tool in the industry of injection molding.

**Author Contributions:** Investigation, J.M.M.-C., D.F.G.-M. and C.M.-D.; resources, J.M.M.-C., D.F.G.-M., B.G.-J. and C.M.-D.; validation, J.M.M.-C., D.F.G.-M., B.G.-J. and C.M.-D.; writing—original draft preparation J.M.M.-C., D.F.G.-M. and C.M.-D.; writing—review and editing, J.M.M.-C., D.F.G.-M. and C.M.-D.; supervision C.M.-D.; project administration C.M.-D.; funding acquisition, C.M.-D. All authors have read and agreed to the published version of the manuscript.

**Funding:** This research work was supported by the University of Jaen with Plan de Apoyo a la Investigación 2021–2022-ACCION1a POAI 2021–2022: TIC-159 and the Plan for Research, Transfer of Knowledge and Entrepreneurship (2022–2023) through the project "Virtual reality as an entrepreneurship tool in industrial design".

**Institutional Review Board Statement:** Not applicable

**Informed Consent Statement:** Not applicable

**Data Availability Statement:** Datasets were analyzed in this study. This data can be found here: [https://youtu.be/SOASYJZ49Nw].

**Acknowledgments:** The authors acknowledge the support of INVELON, who collaborated in the research.

**Conflicts of Interest:** Author Bartolomé Gutierrez-Jiménez was employed by the company Moldes Tuccibérica. The remaining authors declare that the research was conducted in the absence of any commercial or financial relationships that could be construed as a potential conflict of interest.

## References

1.  Torres-Alba, A.; Mercado-Colmenero, J.M.; Caballero-Garcia, J.D.D.; Martin-Doñate, C. Application of New Conformal Cooling Layouts to the Green Injection Molding of Complex Slender Polymeric Parts with High Dimensional Specifications. *Polymers* **2023**, *15*, 558.
2.  Research and Markets. *Injection Molded Plastics Market Size, Share & Trends Analysis Report by Raw Material (Polypropylene, ABS, HDPE, Poly-Styrene), by Application (Packaging, Automotive & Transportation, Medical), by Region, and Segment Forecasts, 2021–2028*; Research and Markets: Dublin, Ireland, year.
3.  Wang, M.L.; Zheng, L.J.; Kang, H.W. 3-Dimenional conformal cooling channel design: Origami-inspired topology optimization approach. *Appl. Therm. Eng.* **2024**, *242*, 122526.
4.  Diaz-Perete, D.; Mercado-Colmenero, J.M.; Valderrama-Zafra, J.M.; Martin-Doñate, C. New procedure for BIM characterization of architectural models manufactured using fused deposition modeling and plastic materials in 4.0 advanced construction environments. *Polymers* **2020**, *12*, 1498.
5.  Arman, S.; Lazoglu, I. A comprehensive review of injection mold cooling by using conformal cooling channels and thermally enhanced molds. *Int. J. Adv. Manuf. Technol.* **2023**, *127*, 2035–2106.
6.  Masato, D.; Kim, S.K. Global Workforce Challenges for the Mold Making and Engineering Industry. *Sustainability* **2023**, *16*, 346.
7.  Wang, X.; Liu, J.; Zhang, Y.; Kristiansen, P.M.; Islam, A.; Gilchrist, M.; Zhang, N. Advances in precision microfabrication through digital light processing: System development, material and applications. *Virtual Phys. Prototyp.* **2023**, *18*, e2248101.
8.  Peixoto, C.; Valentim, P.T.; Sousa, P.C.; Dias, D.; Araújo, C.; Pereira, D.; Machado, C.F.; Pontes, A.J.; Santos, H.; Cruz, S. Injection molding of high-precision optical lenses: A review. *Precis. Eng.* **2022**, *76*, 29–51.
9.  Masato, D.; Piccolo, L.; Lucchetta, G.; Sorgato, M. Texturing Technologies for Plastics Injection Molding: A Review. *Micromachines* **2022**, *13*, 1211.
10. Feng, S.; Kamat, A.M.; Pei, Y. Design and fabrication of conformal cooling channels in molds: Review and progress updates. *Int. J. Heat Mass Transf.* **2021**, *171*, 121082.
11. Formentini, G.; Boix Rodríguez, N.; Favi, C. Design for manufacturing and assembly methods in the product development process of mechanical products: A systematic literature review. *Int. J. Adv. Manuf. Technol.* **2022**, *120*, 4307–4334.
12. Alderighi, T.; Malomo, L.; Auzinger, T.; Bickel, B.; Cignoni, P.; Pietroni, N. State of the art in computational mould design. In *Computer Graphics Forum*; Wiley Online Library: Hoboken, NJ, USA, 2022; Volume 41, pp. 435–452.
13. Asad, M.M.; Naz, A.; Churi, P.; Guerrero AJ, M.; Salameh, A.A. Mix method approach of measuring VR as a pedagogical tool to enhance experimental learning: Motivation from literature survey of previous study. *Educ. Res. Int.* **2022**, *2022*, 8262304.
14. Haleem, A.; Javaid, M.; Singh, R.P.; Rab, S.; Suman, R.; Kumar, L.; Khan, I.H. Exploring the potential of 3D scanning in Industry 4.0: An overview. *Int. J. Cogn. Comput. Eng.* **2022**, *3*, 161–171.
15. Sidani, A.; Dinis, F.M.; Duarte, J.; Sanhudo, L.; Calvetti, D.; Baptista, J.S.; Martins, J.P.; Soeiro, A. Recent tools and techniques of BIM-Based Augmented Reality: A systematic review. *J. Build. Eng.* **2021**, *42*, 102500.
16. Demirel, H.O.; Ahmed, S.; Duffy, V.G. Digital human modeling: A review and reappraisal of origins, present, and expected future methods for representing humans computationally. *Int. J. Hum.-Comput. Interact.* **2022**, *38*, 897–937.
17. de Freitas, F.V.; Gomes MV, M.; Winkler, I. Benefits and challenges of virtual-reality-based industrial usability testing and design reviews: A patents landscape and literature review. *Appl. Sci.* **2022**, *12*, 1755.

18. Mathur, A.; Pirron, M.; Zufferey, D. Interactive programming for parametric cad. In *Computer Graphics Forum*; Wiley Online Library: Hoboken, NJ, USA, 2020; Volume 39, pp. 408–425.

19. Fechter, M.; Schleich, B.; Wartzack, S. Comparative evaluation of WIMP and immersive natural finger interaction: A user study on CAD assembly modeling. *Virtual Real.* **2022**, *26*, 143–158.

20. Zhang, Y.; Kwok, T.H. Design and interaction interface using augmented reality for smart manufacturing. *Procedia Manuf.* **2018**, *26*, 1278–1286.

21. Marougkas, A.; Troussas, C.; Krouska, A.; Sgouropoulou, C. Virtual reality in education: A review of learning theories, approaches and methodologies for the last decade. *Electronics* **2023**, *12*, 2832.

22. Bordegoni, M.; Ferrise, F. Designing interaction with consumer products in a multisensory virtual reality environment. *Virtual Phys. Prototyp.* **2013**, *8*, 51–64.

23. Ferrise, F.; Bordegoni, M. Cugini, Interactive virtual prototypes for testing the interaction with new products. *Comput.-Aided Des. Appl.* **2013**, *10*, 515–525. https://doi.org/10.3722/cadaps.2013.515-525.

24. Gebhardt, S.; Pick, S.; Voet, H.; Utsch, J.; Al Khawli, T.; Eppelt, U.; Reinhard, R.; Büscher, C.; Hentschel, B.; Kuhlen, T.W. flapAssist: How the integration of VR and visualization tools fosters the factory planning process. In Proceedings of the 2015 IEEE Virtual Reality (VR), Arles, France, 23–27 March 2015; Höllerer, T., Ed.; IEEE: New York, NY, USA, 2015; pp. 181–182.

25. Sampaio, A.Z.; Gomes, A.R.; Gomes, A.M.; Santos, J.P.; Rosario, D.P. Collaborative maintenance and construction of buildings supported on virtual reality technology. In Proceedings of the 6th Iberian Conference on Information Systems and Technologies (CISTI 2011), IEEE, Chaves, Portugal, 15–18 June 2011; pp. 1–4.

26. Boton, C. Supporting constructability analysis meetings with immersive virtual reality-based collaborative bim 4d simulation. *Autom. Constr.* **2018**, *96*, 1–15.

27. Sun, S.H.; Tsai, L.Z. Development of virtual training platform of injection molding machine based on VR technology. *Int. J. Adv. Manuf. Technol.* **2012**, *63*, 609–620.

28. Chandra Sekaran, S.; Yap, H.J.; Musa, S.N.; Liew, K.E.; Tan, C.H.; Aman, A. The implementation of virtual reality in digital factory—A comprehensive review. *Int. J. Adv. Manuf. Technol.* **2021**, *115*, 1349–1366.

29. Wolfartsberger, J. Analyzing the potential of Virtual Reality for engineering design review. *Autom. Constr.* **2019**, *104*, 27–37.

30. Aurora, V.R Invelon Technologies SL Polígono Activa Park, Carrer Térmens, 3 25191 - Lleida (Lleida) Spain. Available online: https://invelon.com/auroravr/ (accessed on 25 October 2023).

31. Catia V5–6R2021. 10, rue Marcel Dassault Paris Campus Vélizy-Villacoublay, 78140 France . Available online: https://www.3ds.com/es/ (accessed on 25 October 2023).

32. Faul, F.; Erdfelder, E.; Lang, A.-G.; Buchner, A. G*Power 3: A flexible statistical power analysis program for the social, behavioral, and biomedical sciences. *Behav. Res. Methods* **2007**, *39*, 175–191.

33. Faul, F.; Erdfelder, E.; Buchner, A.; Lang, A.-G. Statistical power analyses using G*Power 3.1: Tests for correlation and regression analyses. *Behav. Res. Methods* **2009**, *41*, 1149–1160.

34. SPSS Statistics Version 27 New Orchard Road Armonk, New York 10504-1722 . Avaible online: https://www.ibm.com/es-es/products/spss-statistics (accessed on 25 October 2023).